\documentclass[11pt]{amsart}

\usepackage{amsmath, amsthm}
\usepackage{amssymb}
\usepackage{amscd}
\usepackage{palatino}
\usepackage{latexsym}
\usepackage{epsfig}
\usepackage{graphics}
\usepackage{amsfonts}
\usepackage{psfrag}
\usepackage{graphicx}

\input xy
\xyoption{all}
\UseComputerModernTips

\usepackage[margin=1in]{geometry}
\numberwithin{equation}{section}
\numberwithin{figure}{section}

\newtheorem{theorem}{Theorem}[section]
\newtheorem{lemma}[theorem]{Lemma}
\newtheorem{proposition}[theorem]{Proposition}

\theoremstyle{definition}
\newtheorem{definition}[theorem]{Definition}

\theoremstyle{remark}

\newtheorem{remark}[theorem]{Remark}

\newcommand{\C}{{\mathbb{C}}}

\newcommand{\Z}{{\mathbb{Z}}}
\newcommand{\Q}{{\mathbb{Q}}}
\newcommand{\R}{{\mathbb{R}}}

\renewcommand{\P}{{\mathbb{P}}}

\renewcommand{\t}{\mathfrak{t}}

\newcommand{\into}{\hookrightarrow}
\newcommand{\onto}{\twoheadrightarrow}

\renewcommand{\mod}{/\!/}

\renewcommand{\emptyset}{\varnothing}

\DeclareMathOperator{\Fred}{Fred}
\DeclareMathOperator{\grad}{grad}
\DeclareMathOperator{\Stab}{Stab}

\DeclareMathOperator{\Hom}{Hom}

\DeclareMathOperator{\codim}{codim}

\DeclareMathOperator{\Crit}{Crit}

\DeclareMathOperator{\pt}{pt}

\DeclareMathOperator{\Span}{span}

\newcommand{\Kernel}{{\mathcal K}}

\begin{document}


\title{The $K$-theory of abelian symplectic quotients}

\author{Megumi Harada}
\address{Department of Mathematics and Statistics, McMaster University, 1280 Main Street West, Hamilton, Ontario L8S 4K1, Canada}
\email{Megumi.Harada@math.mcmaster.ca} 
\urladdr{http://www.math.mcmaster.ca/{}Megumi.Harada/}

\author{Gregory D. Landweber}
\address{Department of Mathematics, Bard College, Annandale-on-Hudson, New York 12504, 
U.S.A.}
\email{landweber@bard.edu}


\date{\today}


\begin{abstract}

  Let $T$ be a compact torus and $(M,\omega)$ a
  Hamiltonian $T$-space.  In a previous paper, the authors showed that
  the $T$-equivariant $K$-theory of the manifold $M$ surjects onto the
  ordinary integral $K$-theory of the symplectic quotient $M \mod T$,
  under certain technical conditions on the moment map. In this paper,
  we use equivariant Morse theory to give a method for computing the
  $K$-theory of $M \mod T$ by obtaining an explicit description of the
  kernel of the surjection \(\kappa: K^*_T(M) \onto K^*(M\mod T).\) 
  Our results
  are $K$-theoretic analogues of the work of Tolman and Weitsman for 
  Borel equivariant cohomology. Further, we prove that under suitable
  technical conditions on the $T$-orbit stratification of $M$, there
  is an explicit Goresky-Kottwitz-MacPherson (``GKM'') type
  combinatorial description of the $K$-theory of a Hamiltonian
  $T$-space in terms of fixed point data. Finally, we illustrate our
  methods by computing the ordinary $K$-theory of compact symplectic
  toric manifolds, which arise as symplectic quotients of an affine
  space $\C^N$ by a linear torus action.

\end{abstract}

\maketitle

\tableofcontents    

\section{Introduction} 

The main result of this manuscript is an explicit description of the
ordinary $K$-theory of an abelian symplectic quotient $M \mod T$ 
in terms of the equivariant
$K$-theory $K^*_T(M)$ of the original Hamiltonian $T$-space $M$, 
where $T$ is a compact torus. 
Symplectic quotients by Lie groups arise naturally in many different fields;
well-known examples are toric varieties and moduli spaces of bundles
over Riemann surfaces. Further, symplectic quotients can often be
identified with Geometric Invariant Theory (``GIT'') quotients in
complex algebraic geometry.  Hence, many moduli spaces that arise as
GIT quotients also have symplectic realizations, and the topological
invariants of such moduli spaces give useful constraints on moduli
problems. In addition, the theory of geometric quantization provides a
fundamental link between the topology of symplectic quotients and
representation theory (see e.g. \cite[Section 7]{Kir98}).  Our
methods, as developed in this manuscript, give a general procedure for
computing the $K$-theory of such spaces when the Lie group is a compact torus. 
By $K$-theory, we mean topological, 
integral $K$-theory, taking $K^0(X)$ to be the 
isomorphism classes of virtual complex vector bundles over $X$ when $X$ is compact, or 
for more general $X$, taking $[X,\Fred(\mathcal{H})]$ for a
complex separable Hilbert space $\mathcal{H}$. 
In the
equivariant case, by $K_T(X)$ we mean Atiyah-Segal $T$-equivariant $K$-theory
\cite{Seg68}, built from $T$-equivariant vector bundles if $X$ is a compact $T$-space, and 
$T$-equivariant maps 
$[X, \Fred(\mathcal{H}_T)]_T$ if $X$ is noncompact (here $\mathcal{H}_T$ contains
every irreducible representation of $T$ with infinite multiplicity, see e.g. \cite{AtiSeg04}).

In the setting of rational Borel equivariant cohomology, a fundamental
symplecto-geometric result of Kirwan \cite{Kir84} states that there is a natural
surjective ring homomorphism
\begin{equation}\label{eq:Kirwan-original}
\kappa_H: H^*_T(M;\Q) \onto H^*(M \mod T; \Q),
\end{equation}
where $H^*_T(M;\Q)$ is the $T$-equivariant cohomology ring of the
original Hamiltonian $T$-space $M$ from which $M\mod T$ is
constructed. Hence, in
order to compute $H^*(M\mod T;\Q)$, it suffices to compute
two objects: the equivariant cohomology ring $H^*_T(M;\Q)$ and the
kernel of $\kappa_H$. For both of these computations, one 
can use {\em equivariant} techniques that are unavailable on the
quotient. 
Following Kirwan's original theorem \cite{Kir84}, this ``Kirwan
method'' has been well developed to yield explicit methods to compute 
both $H^*_T(M;\Q)$ and $\ker(\kappa_{H})$. In this paper, we generalize to $K$-theory the explicit
computation of the kernel of $\kappa_H$ given by Tolman and Weitsman
\cite{TW03}.
Further, we also take to the $K$-theory setting the GKM-type combinatorial
description of the $T$-equivariant cohomology of Hamiltonian
$T$-spaces, which is motivated by the original
work of Goresky, Kottwitz, and MacPherson
\cite{GKM}.

We now explain the setting of our results. 
Suppose that $M$ is a symplectic manifold
with a Hamiltonian $T$-action, i.e., there exists a moment map \(\mu: M
\to \t^*.\) Assuming that $T$ acts freely on the level set
$\mu^{-1}(0)$, the 
symplectic quotient is then defined as \(M \mod T :=
\mu^{-1}(0)/T.\) We wish to compute the ordinary $K$-theory $K^*(M\mod T)$. The first step in
this direction is the $K$-theoretic analogue of the Kirwan
surjectivity result~\eqref{eq:Kirwan-original} above, proven in
\cite{HarLan05}: just as in the cohomology case, there is a natural
ring homomorphism $\kappa$ induced by the natural inclusion
\(\mu^{-1}(0) \into M\) as follows:
\begin{equation}\label{eq:Kirwan-Ktheory}
\xymatrix{
K_{T}^{*}(M) \ar[r] \ar[dr]_{\kappa} & K_{T}^{*}\left(\mu^{-1}(0)\right) \ar[d]^{\cong} \\
& K^{*}(M\mod T),
} 
\end{equation}
and $\kappa$ is surjective. 
The main result of this paper is to give an explicit
computation of the kernel of the $K$-theoretic Kirwan map $\kappa$ above.
This, together with the surjectivity of
$\kappa$, gives us a method to explicitly describe the
$K$-theory of abelian symplectic quotients.

A key element in our arguments is
the $K$-theoretic Atiyah-Bott lemma, which is a $K$-theoretic analogue of a 
fact originally proven in \cite{AtiBot83} in
the setting of Borel equivariant cohomology. In this manuscript,
we use the formulation 
given in \cite[Lemma 2.1]{HarLan05}; a version in the algebraic
category is given in \cite[Lemma 4.2]{VeVi03}. As we noted in \cite{HarLan05},
the Atiyah-Bott lemma is a crucial step in many Morse-theoretic proofs
in symplectic geometry using the moment map $\mu$, and once we have a
$K$-theoretic Atiyah-Bott lemma, it may be expected that many
symplectic-geometric results in the setting of rational
Borel equivariant cohomology carry over to that of (integral)
$K$-theory. 
Indeed, the results of this manuscript can be viewed as 
illustrations of this principle. 

We now briefly outline the contents of this manuscript. 
In Section~\ref{sec:injectivity} we prove an injectivity result, Theorem~\ref{theorem:injectivity},
which states that the $T$-equivariant $K$-theory of the Hamiltonian
$T$-space, under some technical conditions on the moment map, injects
into the $T$-equivariant $K$-theory of the fixed point set $M^T$ via the natural restriction map
\(\imath^*: K^*_T(M) \to K^*_T(M^T).\) 
Similar results were obtained in the compact Hamiltonian setting
already in \cite{GK03}, and in the algebraic category (in algebraic
$K$-theory, for actions
of diagonalizable group schemes on smooth proper schemes over perfect
fields) in
\cite{VeVi03}. Our contribution is to extend the result in the
Hamiltonian setting to a situation where $M$ may not be compact; the
proof uses standard Morse-theoretic techniques in symplectic geometry
using a generic component of the moment map and the $K$-theoretic
Atiyah-Bott lemma mentioned above.  
This injectivity result is an important preliminary step for our main
result 
because the $T$-action on $M^T$ is
trivial, and the ring isomorphism $K^*_T(M^T) \cong K^*(M^T) \otimes K^*_T(\pt)$
makes the equivariant $K$-theory of the fixed point set straightforward to compute. 

With these tools in hand, we prove the following
in Section~\ref{sec:kernel}:

\begin{theorem}\label{theorem:intro-1}
  Let $T$ be a compact torus and $(M,\omega)$ a Hamiltonian $T$-space
  with moment map \(\mu: M \to \t^*.\) Suppose there exists a
  component of the moment map which is proper and bounded below, and
  further suppose that \(M^T\) has only finitely many connected
  components. Let
\begin{equation}\label{eq:def-Z-intro}
Z := \bigl\{\mu(C) \bigm \vert  C \text{ a
  connected component of } \Crit(\|\mu\|^2) \subseteq M \bigr\}
\subseteq \t^* \cong \t
\end{equation} 
be the set of images under $\mu$ of
connected components of the critical set of the norm-square of $\mu$. 
Suppose that $T$ acts freely on $\mu^{-1}(0)$ and let $M \mod T :=
\mu^{-1}(0)/T$ be the symplectic quotient. 
For \(\xi \in \t,\) define
\begin{align*}
M_{\xi} &  :=  \bigl\{ x \in M \bigm\vert \left<\mu(x), \xi\right> \leq 0 \big\}, \\
\Kernel_{\xi} & :=  \bigl\{ \alpha \in K^*_T(M) \bigm\vert  \alpha |_{M_{\xi}} = 0 \bigr\}, \quad \mbox{and} \\
\Kernel & :=  \sum_{\xi \in Z \subseteq \t} \Kernel_{\xi}.
\end{align*}
Then there is a short exact sequence 
\[
\xymatrix{
0 \ar[r] & \Kernel \ar[r] & K_T^*(M) \ar[r]^-{\kappa} & K^*(M \mod T) \ar[r] & 0, 
} 
\]
where \(\kappa: K^*_T(M) \to K^*(M \mod T)\) is the $K$-theoretic Kirwan map. 
\end{theorem}

This theorem is a $K$-theoretic analogue and a slight refinement of the
computation of the kernel of the Kirwan map for rational
Borel equivariant cohomology given by Tolman and Weitsman
\cite{TW03}. In their paper, Tolman and Weitsman give an expression for
the kernel $\ker(\kappa_H)$ which is a sum of ideals
$\Kernel_{\xi}$ (defined similarly to those given above) over {\em all} \(\xi \in \t,\) and it is not
immediately evident that such an infinite sum, in the case when $M$ is
not compact, yields a finite algorithm for the computation of the
kernel. In our version of the computation,
we refine the statement by explicitly exhibiting 
$\Kernel$ as a sum over 
a {\em finite} set $Z$. For a discussion of a different
simplification of the kernel computation for 
rational Borel equivariant cohomology in the compact Hamiltonian case, see \cite{Gol02}.

In Section~\ref{sec:GKM} we prove a $K$-theoretic version of a
Goresky-Kottwitz-MacPherson (``GKM'') combinatorial description
of the $T$-equivariant $K$-theory of a Hamiltonian $T$-space
satisfying certain conditions on the orbit type stratification.
In this setting, we 
give an explicit combinatorial description of the image of $\imath^*$
in \(K^*_T(M^T) \cong K^*(M^T) \otimes K^*_T(\pt).\) 
This result
should be  interpreted as a part of  a large body of work motivated by,
and in many ways generalizing, the original work of Goresky, Kottwitz,
and MacPherson \cite{GKM}, who considered algebraic  torus actions on
projective algebraic  varieties and Borel equivariant  cohomology with
$\C$ coefficients.   For instance,  it is now  known that  similar GKM
results  hold  in  the  setting of  Borel equivariant  cohomology  for
Hamiltonian   $T$-spaces (see e.g. \cite{TW99, HarHol05}).  Similar 
GKM-type 
results for Hamiltonian $T$-spaces in equivariant $K$-theory with $\C$
coefficients are
discussed in \cite{KnuRos03, GK03}, for other equivariant cohomology
theories (under suitable hypotheses on the cohomology theory) and for 
more general $T$-spaces in
\cite{HHH05}, and for equivariant $K$-theory in the algebraic category in \cite{VeVi03}. 
Our  contribution in this  section is to prove, using results in \cite{HHH05}, that such a 
GKM theorem also holds in equivariant $K$-theory (over $\Z$) for Hamiltonian $T$-spaces
which satisfy certain conditions  on the $T$-orbit stratification.

Finally, in Section~\ref{sec:example} we use our methods to give a computation
of the ordinary $K$-theory of smooth compact projective toric
varieties, which can be obtained by a symplectic quotient of an affine
space $\C^N$. 
This rederives, using symplectic
Morse-theoretic techniques, a description of the $K$-theory of these toric
varieties analogous to the Stanley-Reisner presentation
(see e.g. \cite[Section
6.2]{VeVi03}, \cite{BorHor06}). 

The work in this manuscript opens many avenues for future research, of
which we now mention a few examples. First, it would be of interest to
give explicit computations of the $K$-theory of more examples of abelian
symplectic quotients, such as polygon spaces
\cite{HauKnu97, HauKnu98} or, more generally, weight varieties \cite{Gol01}.
Furthermore, although we restrict our attention in this manuscript to
the case where the symplectic quotient $M \mod G$ is a manifold, we
expect that an orbifold version of our results will still hold in the
situation where $G$ acts only locally freely on $\mu^{-1}(0)$, making
the quotient ${\mathcal X} = \mu^{-1}(0)/G$ an orbifold. Here we use
the definition of the ``full orbifold $K$-theory''
\(\mathsf{K}_{\mathrm{orb}}(\cdot)\) of an orbifold
given in \cite{JarKauKim05}, where it is also shown that there is an
orbifold Chern character map from the full orbifold $K$-theory
$\mathsf{K}_{\mathrm{orb}}({\mathcal X})$ to the Chen-Ruan orbifold
cohomology $H^*_{CR}({\mathcal X})$ of the orbifold ${\mathcal
X}$.
Methods for computing the Chen-Ruan orbifold cohomology of
orbifold symplectic quotients were given in \cite{GHK05} by using
Kirwan surjectivity methods in addition to explicit computations of
the kernel of the Kirwan map in the manifold case.  We expect that,
using the results of this manuscript and an approach similar to
\cite{GHK05}, we can also compute the full orbifold $K$-theory of
orbifold symplectic quotients. Finally, it would be of interest to
prove a $K$-theoretic analogue of the simplification of the
computation of $\ker(\kappa)$ given by Goldin in \cite{Gol02} in the
case of rational Borel equivariant cohomology. Goldin restricted her
considerations to the case of {\em compact} Hamiltonian spaces in
\cite{Gol02}, but we expect that a similar statement should still
hold.  We intend to explore these and related topics in future work.

\bigskip

\noindent {\bf Acknowledgements.}  We thank Jonathan Weitsman for helpful discussions. 
The second author thanks the University of
Toronto and the Fields Institute for their hospitality and support
while conducting a portion of this research. Both authors thank the
American Institute of Mathematics and the 
Banff International Research Station for their 
hospitality. 

\section{Injectivity into the $K$-theory of the fixed points}\label{sec:injectivity}

We begin with a $K$-theoretic version of the injectivity theorem of
Kirwan, which is a key technical tool that we will need in the later sections. 
Suppose given $M$ a Hamiltonian $T$-space with moment map
\(\mu: M \to \t^*,\) with a component which is proper and bounded
below. In this section only, we will for notational convenience use
$\Sigma$ to denote the $T$-fixed point set $M^T$ of $M$. We
additionally assume that $\Sigma$ has only finitely many connected
components. In this situation, it is well-known that 
the inclusion \(\imath: \Sigma \into M\)
induces a map in rational Borel equivariant cohomology 
\begin{equation}\label{eq:injectivity-cohomology} 
\xymatrix{
\imath^*: H^*_T(M;\Q) \ar[r] &  H^*_T(\Sigma;\Q),
}
\end{equation}
which is an injection. We will prove a $K$-theoretic version of this
injectivity~\eqref{eq:injectivity} by Morse theory using a generic
component of the $T$-moment map $\mu$.  Our proof will follow that
given in \cite[Theorem 2.6]{HarHol05}.

Recall that the 
components $\mu^\xi := \left< \mu, \xi\right>$ for $\xi\in\t$ of the moment map are
Morse-Bott functions on $M$. We call a component {\em generic} if the
critical set of $\mu^\xi$ is precisely the fixed point set $\Sigma$. This
is true of an open dense set of directions \(\xi \in \t.\) 
Our assumption above that there exists a component which is proper and
bounded below guarantees that our Morse-theoretic arguments
will work.

In the course of the proof we will use a special case of the $K$-theoretic Atiyah-Bott lemma,
which we restate here for reference.

\begin{lemma}\label{lemma:AtiyahBott} (\cite[Lemma
  4.2]{VeVi03}, \cite[Lemma 2.3]{HarLan05})
Let a compact connected Lie group $G$ act fiberwise linearly on a
complex vector bundle \(\pi: E \to X\) over a compact connected
$G$-manifold $X$. Assume that a circle subgroup \(S^1 \subseteq G\)
acts on $E$ so that the fixed point set is precisely the zero section
$X$. Choose an invariant metric on $E$ and let $D(E)$ and $S(E)$ denote the
disc and sphere bundles, respectively. Then the long exact sequence
for the pair \((D(E),S(E))\) in equivariant $K$-theory splits into short
exact sequences
\[
\xymatrix{ 
0 \ar[r] & K^*_G\bigl(D(E),S(E)\bigr) \ar[r] & K^*_G\bigl(D(E)\bigr) \ar[r] & K^*_G\bigl(S(E)\bigr) \ar[r] & 0.
}
\]
\end{lemma} 

\begin{remark} 
  Recall that the splitting of the long exact sequence is equivalent,
  by the Thom isomorphism, to the statement that the $K$-theoretic
  equivariant Euler class of $E$ is not a zero divisor in $K^*_G(X)$.
  We will mainly use this latter perspective.
\end{remark}

In what follows, we need only consider the Atiyah-Bott lemma for
the case \(G = T.\) Moreover, in our applications of
Lemma~\ref{lemma:AtiyahBott}, the $T$-action on the base manifold $X$
will be trivial (since the base manifold is a component of a $T$-fixed
set).  Hence, in order to use the lemma, we need only check
that there exists a subgroup $S^1 \subseteq T$ acting on $E$ fixing
precisely the zero section.

Before proceeding to the main argument, we prove the following
simple technical lemma:

\begin{lemma}\label{lemma:genericComponent}
  Let $T$ be a compact torus and $(M,\omega)$ a Hamiltonian $T$-space
  with moment map $\mu: M \to \t^*$. Suppose there exists a component
  of $\mu$ which is proper and bounded below, and further suppose that
  $\Sigma$ has only finitely many connected components. Then
  there exists a generic component \(f := \mu^{\xi}\) of $\mu$ which
  is proper, bounded below, and with respect to which all components
  of $\Sigma$ have different critical values.
\end{lemma}

\begin{proof}
  If there is a component $\mu^{\eta}$ of $\mu$ which is proper and
  bounded below, then $\mu$ itself is proper, since a component is
  given by a linear projection.  A result of Lerman, Meinrenken,
  Tolman, and Woodward \cite[Theorem 4.2]{LMTW} states that in this
  situation the image $\mu(M)$ of $M$ is convex and locally
  polyhedral.  Thus, by taking a small enough perturbation $\xi$ of $\eta$,
we may arrange that all 
  components of the fixed point set $M^T$ have different critical
  values with respect to $\mu^{\xi}$ and still have $\mu^{\xi}$ proper
  and bounded below. 
\end{proof}

With this lemma in hand, we may now state the proposition
which is our main technical tool to
show injectivity. This is the $K$-theoretic analogue of
\cite[Proposition~2.4]{HarHol05}.

\begin{proposition}\label{prop:Morse}
Let $T$ be a compact torus and $(M,\omega)$ a Hamiltonian $T$-space
with moment map \(\mu: M \to \t^*.\) Suppose there exists 
a component of the moment map
which is proper and bounded below, and further suppose that $\Sigma$
has only finitely many connected components. Let $f := \mu^{\xi}$ be a
generic component of $\mu$ as in
Lemma~\ref{lemma:genericComponent}. 
Let $c$ be a critical value of $f$, and pick $\varepsilon >0$ such that
$c$ is the only critical value of $f$ in $(c-\varepsilon,c+\varepsilon)$. Let
$\Sigma_c$ be the component of $\Sigma$ with
\(f(\Sigma_c) = c,\) and let 
\(M_c^{\pm} := f^{-1}(-\infty, c \pm \varepsilon).\) Then the long
exact sequence of the pair \((M_c^+, M_c^-)\) splits into short exact
sequences
\[
\xymatrix{
0 \ar[r] & K^*_T\bigl(M_c^+, M_c^-\bigr) \ar[r] & K^*_T\bigl(M_c^+\bigr) \ar[r]^-{k^*} & 
K^*_T\bigl(M_c^-\bigr) \ar[r] & 0,
}
\]
where \(k: M_c^- \into M_c^+\) is the inclusion. Moreover, the restriction map
\(K^*_T(M_c^+) \to K^*_T(\Sigma_c),\) given by the inclusion
\(\Sigma_c \into M_c^+,\) induces an isomorphism from the kernel of
$k^*$ to the classes of $K^*_T(\Sigma_c)$ that are multiples of the
equivariant Euler class $e(\Sigma_c)$ of the negative normal bundle to
$\Sigma_c$ with respect to $f$.
\end{proposition}

\begin{proof} 
This is a $K$-theoretic version of the proof given in \cite{TW99}.
Let $D(\nu_c)$ and $S(\nu_c)$ denote the disc and sphere bundles of the negative
normal bundle $\nu_c$ to the fixed set $\Sigma_c$ with respect to
$f$. By the assumption on $\varepsilon$, the component $\Sigma_c$ is the only
component of $\Crit(f)$ in $f^{-1}(c-\varepsilon,
c+\varepsilon)$. Similarly since $f$ is proper and $T$-invariant, the
negative gradient flow with respect to $f$ using a $T$-invariant
metric gives a $T$-equivariant
retraction of the pair \((M_c^+, M_c^-)\) to the pair $(D(\nu_c),
S(\nu_c))$. Using 
the $T$-equivariant Thom isomorphism, we conclude that the long exact
sequence of the pair 
\[
\cdots \to K^*_T(M^+_c, M^-_c) \to K^*_T(M_c^+) \to K^*_T(M_c^-) \to
\cdots
\]
splits into short exact sequences if and only if the equivariant Euler
class $e(\nu_c)$ is not a zero divisor. (We also note that, unlike the
rational cohomology case, since 
$f$ is a component of a moment map, the Morse index of $f$ is
even. Since 
\(K^{*-\lambda}_T(D(\nu_c)) \cong K^*_T(D(\nu_c))\) by Bott
periodicity, there is no degree shift here.) Moreover, the negative
normal bundle $E$ is complex by a result of Kirwan \cite{Kir84}. The
group $T$ fixes exactly the zero section $\Sigma_c$ because $\Sigma_c$
is defined to be a component of the $T$-fixed point set. We may apply
Lemma~\ref{lemma:AtiyahBott} taking \(G=T, X=\Sigma_c,\) and \(E=\nu_c,\) together with a 
suitable choice of $S^1 \subset T$, to conclude that $e(\Sigma_c)$ is not
a zero divisor. The proposition follows.
\end{proof}

We now come to the main theorem of this section. 
With Lemmas~\ref{lemma:AtiyahBott} and~\ref{lemma:genericComponent}
and Proposition~\ref{prop:Morse} in hand, the proof now exactly
follows that given for \cite[Theorem 2.5]{HarHol05}, so we will not
reproduce it here. 

\begin{theorem}\label{theorem:injectivity} 
  Let $T$ be a compact torus and $(M,\omega)$ a Hamiltonian $T$-space
  with moment map \(\mu: M \to \t^*.\) Suppose there exists a
  component of the moment map which is proper and bounded below, and
  further suppose that $\Sigma$ has only finitely many connected
  components. Let $f := \mu^{\xi}$ be a generic component of $\mu$ as
  in Lemma~\ref{lemma:genericComponent}.  Let \(\imath: \Sigma \into
  M\) be the inclusion of the fixed point set into $M$. Then the
  restriction map
\begin{equation}\label{eq:injectivity}
\imath^*: K^*_T(M) \to K^*_T(\Sigma) 
\end{equation}
is injective. 
\end{theorem}

\begin{remark}
Since the Atiyah-Bott lemma and Proposition~\ref{prop:Morse} tell us
that the Euler classes $e(\Sigma_c)$ of the negative normal bundles to the
fixed point components are not zero divisors, in order to obtain the
injectivity~\eqref{eq:injectivity}, we could also take the Morse
stratification \(M = \bigsqcup_c S_c,\) where the $S_c$ are the flow-down
manifolds associated to $\Sigma_c$, and apply the injectivity theorem
\cite[Theorem 2.3]{HHH05}, which holds for more general equivariant
cohomology theories. We present this straightforward Morse-theoretic 
argument
here since we expect it to be more familiar to some readers.
\end{remark}

\section{The kernel of the Kirwan map}\label{sec:kernel}

We have already observed in the introduction that 
in order to make explicit computations of the $K$-theory rings of abelian
symplectic quotients, it is necessary to explicitly identify the kernel of the
Kirwan map \(\kappa: K_T^*(M) \to K^*(M \mod T).\) 
In the setting of rational Borel equivariant
cohomology, Tolman and Weitsman showed in \cite{TW03} that the kernel of
$\kappa$ can be identified as a sum of certain ideals in
$H^*_T(M;\Q)$. 
In this section we state and prove the (integral) $K$-theoretic analogue of
this result. This will allow us to make the 
explicit computation in Section~\ref{sec:example}.

Before proceeding to the main
theorem, we take a moment to discuss the relationships between
the various technical Morse-theoretic hypotheses on a moment map $\mu$ which are
used in the literature. These are:
\begin{enumerate}
\item there exists a component $\mu^{\xi}$ of the moment map which is proper
and bounded below, 
\item $\mu$ is proper, and
\item \(\|\mu\|^2\) is proper.
\end{enumerate}
It is straightforward to see that (1) implies (2), and that
(2) and (3) are equivalent. When the Hamiltonian space $M$ is compact,
then any of these conditions automatically holds; the point is that
when $M$ is not compact, these conditions ensure that the
Morse-theoretic arguments given below using $\mu$ still work. 
For the proof of the $K$-theoretic Kirwan surjectivity
theorem in \cite{HarLan05}, we need only the second (or equivalently
third) condition, namely, that $\mu$ is proper. However, for our explicit kernel computations 
we will be
using Theorem~\ref{theorem:injectivity} which requires the 
hypothesis (1) and additionally the condition that $M^T$ has finitely
many components, so we continue
making these stronger assumptions. In practice this is not very
restrictive.\footnote{For instance, although hyperk\"ahler quotients are rarely
compact, they often admit natural Hamiltonian torus actions which do
satisfy (1); for example, for hypertoric varieties, this fact was
exploited in \cite{HarHol05} to give combinatorial descriptions of
their equivariant cohomology rings, and for hyperpolygon spaces, Konno uses a
Hamiltonian $S^1$-action satisfying (1) in his proof of Kirwan surjectivity
for these spaces \cite{Kon02}. (See also \cite{HauPro05}
for more on examples of this type.)} Moreover, our hypotheses guarantee that the
algorithm we present for computing the kernel of $\kappa$ is {\em
finite}, as will be seen in Lemma~\ref{lemma:normsquareCrit}.

We first give the statement of the main theorem of this section.
In the following and in the sequel, we identify \(\t \cong \t^*\)
using an inner product on $\t$, so the moment map $\mu$
can be considered as a map taking values in $\t$. 

\begin{theorem}\label{theorem:kernelforT}
  Let $T$ be a compact torus and $(M,\omega)$ a Hamiltonian $T$-space
  with moment map \(\mu: M \to \t^*.\) Suppose there exists a
  component of the moment map which is proper and bounded below, and
  further suppose that \(M^T\) has only finitely many
  connected components. Let
\begin{equation}\label{eq:Z-def}
Z := \bigl\{\mu(C) \bigm\vert  C \text{ a
  connected component of }  \Crit(\|\mu\|^2) \subseteq M\bigr\}
\subseteq \t^* \cong \t
\end{equation}
be the set of images under $\mu$ of
components of the critical set of $\|\mu\|^2$. 
Suppose that $T$ acts freely on the level set $\mu^{-1}(0)$, and let
$M \mod T$ be the symplectic quotient. 
 For \(\xi \in \t,\) define
\begin{align*}\label{eq:Mxi}
M_{\xi} &  :=  \bigl\{ x \in M \bigm\vert \left<\mu(x), \xi\right> \leq 0 \bigr\}, \\
\Kernel_{\xi} & :=  \bigl\{ \alpha \in K^*_T(M) \bigm\vert  \alpha |_{M_{\xi}} = 0 \bigr\}, \quad \mbox{and} \\
\Kernel & :=  \sum_{\xi \in Z} \Kernel_{\xi}.
\end{align*}
Then there is a short exact sequence 
\[
\xymatrix{
0 \ar[r] & \Kernel \ar[r] & K_T^*(M) \ar[r]^-{\kappa} & K^*(M \mod T) \ar[r] & 0, 
} 
\]
where \(\kappa: K^*_T(M) \to K^*(M \mod T)\) is the Kirwan map. 
\end{theorem}

We note that Theorem~\ref{theorem:kernelforT} is both a $K$-theoretic
analogue of the kernel computation in $H^*_T(-;\Q)$ by Tolman and
Weitsman \cite[Theorem 4]{TW03} and a slight refinement of their
statement, as we now explain. In \cite[Theorem 4]{TW03}, the
kernel is expressed as a sum over the infinite set of all elements in
$\t$. In the case where $M$ is compact, they comment that their
algorithm is in fact finitely computable in \cite[Remark
5.3]{TW03}, but they do not explicitly address what happens in the
non-compact situation. In Theorem~\ref{theorem:kernelforT} we have
expressed $\Kernel$ as a sum over a certain set $Z$, which is shown to
be finite in Lemma~\ref{lemma:normsquareCrit} below. The point is
that our $K$-theoretic kernel $\Kernel$ is explicitly computed
by a finite algorithm even in the non-compact case. 

We now show that $Z$ is finite. Recall that under our properness
hypotheses, both $\mu$ and $\|\mu\|^2$ are proper. 

\begin{lemma}\label{lemma:normsquareCrit}
Let $T$ be a compact torus and $(M,\omega)$ a Hamiltonian $T$-space
with moment map \(\mu: M \to \t^*.\) Suppose there exists a component
of the moment map which is proper and bounded below, and further
suppose that \(M^T\) has only finitely many connected
components. Then there are only finitely many connected components of
the critical set of the norm-square \(\|\mu\|^2: M \to \R_{\geq 0}\)
of the moment map. 
In particular, the set $Z$ in~\eqref{eq:Z-def} is finite. 
\end{lemma}

\begin{proof} 
For a subgroup $H \subseteq T$ of the torus $T$, let $M_H$ denote the
subset of $M$ consisting of points $p$ with
\(\Stab(p) = H.\) Let 
\begin{equation}\label{eq:orbit-stratification}
M = \bigsqcup_{H \subseteq T} M_H
\end{equation}
be the orbit type stratification of $M$ as in \cite[Section
3.5]{GGK}. Since $T$ is compact, the action is proper. By standard
Hamiltonian geometry \cite{CdS01} and the theory of proper group
actions (see e.g. \cite[Appendix B]{GGK}), the closure
$\overline{M_H^{(i)}}$ of each connected component $M_H^{(i)}$ of
$M_H$ is itself a closed Hamiltonian $T$-space with moment map given
by restriction. Since the restriction of $\Phi$ to
$\overline{M_H^{(i)}}$ is also proper and bounded below, each
$\overline{M_H^{(i)}}$ contains a $T$-fixed point. By equivariant
Darboux, any $T$-invariant tubular
neighborhood of a connected component of $M^T$ intersects only
finitely many $T$-orbit types. Since by assumption there are only
finitely many connected components of $M^T$, we may conclude that
there are only finitely many orbit types in the
decomposition~\eqref{eq:orbit-stratification}, and that each
$M_H$ has only finitely many connected components. 
Kirwan proves \cite[Lemma 3.12]{Kir84} that there is at most one
critical value of \(\|\mu\|^2\) on each \(\Phi((M_H)^{(i)}),\) so  
we may conclude that there are only finitely many critical values of
$\|\mu\|^2$. Since $\mu$ is proper, there are only finitely many
connected components of $\Crit(\|\mu\|^2)$. The final assertion in the
theorem follows from \cite[Corollary 3.16]{Kir84}. 
\end{proof}

The proof of Theorem~\ref{theorem:kernelforT} closely follows the
Morse-theoretic argument given in the compact case in $H^*_T(-;\Q)$ in
\cite[Theorem 3]{TW03}, except that we use our $K$-theoretic
Atiyah-Bott lemma (Lemma~\ref{lemma:AtiyahBott}). Hence we do not give all the details
below. Instead, we only briefly indicate why 
the sum can in fact be restricted to the set $Z$. 

\begin{proof}[Proof of Theorem~\ref{theorem:kernelforT}]

We refer the reader to \cite{TW03} for details. First, the only
point of the proof requiring substantial argument is to prove that
\(\ker(\kappa) \subseteq \Kernel.\) Second, we order the connected
components of the critical sets of the norm-square $\|\mu\|^2$; denote
these as \(\{C_i\}_{i=0}^m,\) where \(C_0 = \mu^{-1}(0).\) Hence
\(\alpha \in \ker(\kappa)\) exactly means \(\alpha |_{C_0} = 0.\)
Third, by observing that all components of $M^T$ are also components
of $\Crit(\|\mu\|^2)$, it suffices by
Theorem~\ref{theorem:injectivity} to show that there exists a \(\beta
\in \Kernel\) such that \(\alpha|_{C_i} = \beta|_{C_i}\) for all
$i$. The final and most important step in this argument is an
inductive construction of the $\beta$, for which it suffices to show
that for \(1 \leq \ell \leq m\) and \(\alpha \in \ker(\kappa)\) such
that \(\alpha|_{C_i} = 0\) for all \(0 \leq i \leq \ell -1,\) there
exists an element \(\alpha' \in \Kernel\) such that \(\alpha'|_{C_i} =
\alpha|_{C_i}\) for all \(0 \leq i \leq \ell.\) The construction given
by Tolman and Weitsman can now be explicitly seen to produce an
element $\alpha$ which is in fact contained in
$\Kernel_{\mu(C_\ell)}$. Since \(\mu(C_{\ell}) \in Z\) by
definition~\eqref{eq:Z-def}, this implies that \(\alpha' \in
\Kernel.\) The rest of the argument follows that in \cite{TW03}. 

\end{proof}

\section{GKM theory in $K$-theory}\label{sec:GKM}

The main result of this section  is to show that, under some technical
conditions  on   the  orbit   stratification,  the  $K$-theory   of  a
Hamiltonian $T$-space has a combinatorial description, using data from
the  equivariant  one-skeleton.  Such  a  description  is  useful  for
explicit computations  of the  kernel of the  Kirwan map.  
As mentioned in the Introduction, this result is best viewed as part
of a large body of work inspired by the original paper of Goresky,
Kottwitz, and MacPherson \cite{GKM}. 
Both the statement of our theorem and its proof 
are $K$-theoretic  versions of those given in \cite[Theorem
  2.11]{HarHol05},  with slight 
differences  in   the  hypotheses
(explained below). 

We begin by defining an extra condition on the $T$-action that will
be necessary to state the theorem.

\begin{definition}\label{def:GKM}
Let $T$ be a compact torus and $(M,\omega)$ a Hamiltonian $T$-space. We
say that the action is {\bf GKM} if $M^T$ consists of finitely many
isolated points, and the $T$-isotropy weights
at each fixed point \(p \in M^T\) are
pairwise linearly independent in $\t^*_{\Z}$.
\end{definition}

\begin{remark}
  The GKM condition on the $T$-isotropy weights given in
  Definition~\ref{def:GKM} is less restrictive than the hypothesis in
  \cite[Theorem~2.11]{HarHol05} that the $T$-weights are relatively
  prime in $H^*_T(\pt;\Z)$. This is due to the difference between
  equivariant Euler classes in $K^*_T(\pt)$ and those in $H^*_T(\pt;\Z)$.
  A more detailed discussion of the differences between the
  Atiyah-Bott lemma in equivariant $K$-theory and in integral Borel equivariant
  cohomology can be found in \cite[Section 2]{HarLan05}.
\end{remark}

We now briefly recall the construction of the combinatorial and
graph-theoretic data used in our theorem. Let $N$ denote the subset of $M$ given by
\[
N := \bigl\{ p \in M \bigm\vert \codim(\Stab(p)) = 1\bigr\}.
\]
Thus $N$ consists of the points in $M$ whose $T$-orbits are exactly
one-dimensional. The {\em equivariant one-skeleton of $M$} is then
defined to be the closure $\overline{N}$ of $N$. Hence
\[
\overline{N} = \bigl\{ p \in M \bigm\vert \codim(\Stab(p)) \leq 1\bigr\} = N \cup M^T.
\]
The GKM condition states that at each fixed point \(p \in M^T,\) the
$T$-isotropy weights are pairwise linearly independent in $\t^*_{\Z}$. Each $T$-weight
space corresponds to a component of $N$, the closure of which is
either an $S^2 \cong \P^1$ (with north and south poles being $T$-fixed
points \(\{p,q\} \subseteq M^T\)) or a copy of $\C$ (with origin a $T$-fixed point). Thus, given
the GKM condition, the one-skeleton $\overline{N}$ is a collection of
projective spaces $\P^1$ and affine spaces $\C$, glued at $T$-fixed
points.
By definition, each component
of $N$ is equipped with a $T$-action which is specified by a weight in
$\t^*_{\Z} \cong \Hom(T,S^1)$; this weight appears in the $T$-weight
decomposition of the isotropy action on the corresponding $T$-fixed
point. (There is a sign ambiguity in the $T$-weight for a component of
$N$ whose closure is a $\P^1$, depending on the choice of north or
south pole; however, this choice does not affect the GKM computation
to be described below.)

From this data we construct the {\em GKM graph}, a labelled graph
$\Gamma = (V,E,\alpha)$, associated to the
GKM $T$-space $M$. The vertices $V$ of $\Gamma$ are the $T$-fixed
points $V = M^T$, and there is an edge $(p,q) \in E$ exactly when
there exists an embedded $\P^1 \subset \overline{N}$ containing as
its two $T$-fixed points $\{p,q\} \subset \P^1$. Additionally, we
label each edge $(p,q)$ with the weight $\alpha_{(p,q)}$ specifying
the $T$-action on the corresponding $\P^1$ as discussed above. 
Note that the components of $N$ corresponding to an affine space $\C$
do not contribute to the GKM graph since each such $\C$ equivariantly retracts to
its corresponding $T$-fixed point.

Now we return to the setting of the previous sections. Let
$(M,\omega)$ be a symplectic manifold equipped with a Hamiltonian
$T$-action. As before, we assume that the moment map \(\mu: M \to \t^*\) has a
component which is proper and bounded below. In addition, we now
assume that the $T$-action on $M$ is GKM, so in particular $M^T$ is a
finite set of isolated points in $M$.
Let $\mu^{\xi}: M \to \R$ be a generic component of $\mu$ which is
proper and bounded below, as chosen in
Lemma~\ref{lemma:genericComponent}. Then $\mu^{\xi}$ is a Morse
function on $M$. Order the critical points $\Crit(\mu^{\xi}) = M^T =
\{p_i\}_{i=1}^m$ so that \(\mu^{\xi}(p_j) < \mu^{\xi}(p_k)\) if and
only if \(j<k.\) Since the Morse function is $T$-invariant,
the negative gradient flow with respect to a $T$-invariant metric is
$T$-equivariant.

Let $U_i$ be the flow-down cell from the critical point $p_i$. Then 
the negative gradient flow gives a $T$-equivariant deformation
retraction of $M$ to 
the union of the $U_i$, i.e.,
\[
M \sim \bigcup_i U_i
\]
is a $T$-equivariant homotopy equivalence. Hence to study the
$T$-equivariant $K$-theory of $M$, we may instead study that of $M' :=
\bigcup_i U_i$.  Each $U_i$ contains a single critical point $p_i$, and
$T_{p_i}U_i$ is a $T$-representation; by the GKM assumption, 
$T$-weights occurring in $T_{p_i}U_i$ are pairwise linearly independent. 

We may now use the results of \cite{HHH05} to show that, under the
conditions outlined above, the GKM graph $\Gamma$ combinatorially
encodes the equivariant $K$-theory of the Hamiltonian $T$-space $M$.
We first state a crucial lemma which involves
$K$-theoretic equivariant Euler classes in $K^*_T(\pt) \cong R(T)$, where $R(T)$ is the representation ring of $T$;
this will be the key step in the proof
of the main theorem.
Let $\C_{\sigma}$ be a $1$-dimensional representation of $T$ with
weight $\sigma \in \t_{\Z}^*$. Recall that the $K$-theoretic
$T$-equivariant Euler class of the $T$-bundle \(\C_{\sigma} \to
\mathrm{pt}\) is 
\[
e_T(\sigma) := 1 - e^{-\sigma} \in K_T^*(\pt) \cong R(T).
\]
By properties of the Euler
class, if \(E = \bigoplus_i \C_{\sigma_i}\) is a direct sum of such
$1$-dimensional representations, then the equivariant Euler class of
$E$ is the product \(e_T(E) = \prod_i (1 - e^{-\sigma_i}).\) The following is a
special case of \cite[Lemma 4.9]{VeVi03}:

\begin{lemma}\label{lemma:relatively_prime}
  Let \(\sigma, \tau\in \t_{\Z}^*\) be linearly independent in
  $\t^*_{\Z}$. Then the corresponding Euler classes \(e_T(\sigma) = 1 -
  e^{-i\sigma}\) and \(e_T(\tau) = 1 - e^{-i\tau} \cong R(T)\) are relatively prime in
  $K^*_T(\pt) = R(T)$.
\end{lemma}

We will use this lemma to obtain the following combinatorial
description of $K^*_T(M)$. Since $M^T$ consists of finitely many
isolated points, 
\[
K^*_T(M^T) \cong \bigoplus_{p \in M^T} K^*_T(p) \cong \left\{h: V = M^T \to
K^*_T(\pt) \cong R(T) \right\}.
\]
We then define the {\em $\Gamma$-subring} of $K^*_T(M^T)$ 
to be
\[
K^*(\Gamma, \alpha) := \left\{ h:V\to K_T^*(\pt) \cong R(T) \ \left| \
\begin{array}{c}h(p)-h(q)\equiv 0\ ({\mathrm{ mod}}\ e(\alpha_{(p,q)})) \\
\mbox{ for every edge } (p,q) \in E
\end{array}\right.
\right\} \subseteq K_T^*(M^T).
\]
We have the following:

\begin{theorem}\label{theorem:GKM}
Let $T$ be a compact torus and $(M,\omega)$ a Hamiltonian $T$-space
with moment map \(\mu: M \to \t^*.\) Suppose there exists a component
of the moment map which is proper and bounded below, and that the $T$-action on
$M$ is GKM.  Then 
the inclusion \(\imath: M^T \into M\) induces an isomorphism
\[
\imath^*: K^*_T(M) \to K^*(\Gamma,\alpha) \subseteq K^*_T(M^T). 
\]
\end{theorem}

\begin{proof} 
  We will use a special case of \cite[Theorem 3.1]{HHH05}, which
  states that for $T$-spaces satisfying certain assumptions, the
  equivariant $K^*_T$-theory is isomorphic via $\imath^*$ to the
  $\Gamma$-subring defined above. 
 We must therefore check that each of
  the hypotheses necessary for this theorem is satisfied.

We have already seen that it suffices to compute the $K^*_T$-theory of
the subspace $M' := \bigcup_i U_i$ of $M$. Let \(M_i := \bigcup_{1 \leq
j \leq i} U_i\) be the union of these cells up to the $i$-th
cell. Then \(M = \bigcup_i M_i\) is a $T$-invariant stratification of
$M'$, and by construction each quotient \(M_i/M_{i-1}\) is homeomorphic to the Thom
space of the $T$-equivariant negative normal bundle \(\nu(p_i) \to
p_i\) to the fixed point $p_i$ with respect to $\mu^{\xi}$.

Each negative normal bundle $\nu(p_i)$ is $T$-orientable, since it is
complex. Moreover, $\nu(p_i)$ decomposes into a direct sum of
one-dimensional $T$-weight spaces $\C_{\alpha_{j,p_i}}$. Identify the
negative normal bundle with \(U_i = M_i \setminus M_{i-1}.\) 
Each $\C_{\alpha_{j, p_i}}$ corresponds to a
component of $N$, so in particular there is an attaching map of the
sphere bundle \(S(\nu(p_i)) \to M_{i-1}\) which, when restricted to
\(S(\C_{\alpha_{j, p_i}}),\) maps it to a critical point $p_j$, where
\(j<i.\) Such attaching maps exactly correspond to the presence of the
embedded $\P^1$s in the equivariant one-skeleton.

Moreover, since the $T$-isotropy weights $\alpha_{j,p_i}$ are all non-zero,
their 
equivariant Euler classes \(1 - e^{-\alpha_{j,p_i}}\) in $K^*_T(\pt)
\cong R(T)$ are not zero divisors in
$K^*_T(\pt)$. Furthermore, by the GKM assumption, the weights at a given critical
point $p_i$ are also pairwise linearly independent. By
Lemma~\ref{lemma:relatively_prime}, this implies that their
equivariant Euler classes are pairwise relatively prime in
$K^*_T(\pt)$.

By \cite[Theorem 3.1]{HHH05}, we may conclude that $K^*_T(M) \cong
K^*_T(M')$ injects via $\imath^*$ into $K^*_T(M^T)$ with image exactly
the $\Gamma$-subring. This concludes the proof. 
\end{proof}

\section{Example: symplectic toric manifolds}\label{sec:example}

In this section, we illustrate the kernel computation procedure
outlined in Section~\ref{sec:kernel} by giving an explicit computation
of the ordinary $K$-theory of compact symplectic toric manifolds. We
thus obtain a description analogous to the Stanley-Reisner
presentation of the (equivariant) cohomology ring of toric varieties
(see e.g. \cite{VeVi03}, \cite{BorHor06}). The Delzant construction \cite{Del88} is
an explicit construction of a compact symplectic toric manifold,
combinatorially specified by a compact Delzant polytope \(\Delta
\subseteq \R^n,\) as a symplectic quotient of an affine space $\C^N$
by a linear (and hence Hamiltonian) $T^k$-action; this allows us to
use Theorem~\ref{theorem:kernelforT}. We only briefly recall the basic
ingredients of this Delzant construction below and refer the reader to
\cite{CdS01} for details.

Let $\Delta \subseteq \R^n$ be a compact Delzant polytope, i.e., it is compact, simple, rational, and
smooth. Let $a_i$ be the primitive outward-pointing normal vectors to
the facets \(\{F_i\}_{i=1}^N\) of $\Delta$. Hence 
\begin{equation}\label{eq:def-Delta}
\Delta = \left\{ x \in \R^n \, \mid \, \langle x, a_i \rangle \leq
\eta_i, \, 1 \leq i \leq N \right\}, 
\end{equation}
for some \(\eta_i \in \R.\) Here \(\langle \, , \, \rangle\) denotes the standard 
inner product on $\R^n$. (As in Section~\ref{sec:kernel}, we identify $\R^n$ with its dual using the inner product.) Consider the torus $T^N$ acting standardly
on $\C^N$. The Lie algebra \(\t^N \cong \R^N\) has standard basis
\(\{\varepsilon_i\}_{i=1}^N\) and its dual \((\t^N)^*\) has dual basis
\(\{u_i\}_{i=1}^N.\) Let \(\beta: \t^N \to \t^n \cong \R^n\) be the
linear map defined by \(\beta(\varepsilon_i) = a_i.\) This gives us an
exact sequence of Lie algebras 
\begin{equation}\label{eq:Delzant-Lie}
\xymatrix{
0 \ar[r] &  \t^k \ar[r]^{\iota} &  \t^N \ar[r]^{\beta} & \t^n \ar[r] &  0,
}
\end{equation}
where $\iota$ is the inclusion of \(\t^k := \ker(\beta)\) into $\t^N$,
with \(k = N - n.\) There is a corresponding exact sequence of Lie
groups
\[ 
1 \to T^k \to T^N \to T^n \to 1.
\]
Hence $\C^N$ is also a Hamiltonian $T^k$-space, given by restricting the action to the subgroup $T^k$. A moment map for this $T^k$-action is given by 
\[
\Phi(z_1, \ldots, z_n) = - \left( \frac{1}{2} \sum_{i=1}^N
|z_i|^2 \iota^* u_i \right) + \iota^* \eta, 
\]
where \(\eta = (\eta_1, \ldots, \eta_N) \in \R^N \cong (\t^N)^*\) is
constructed from the constants $\eta_i$ in~\eqref{eq:def-Delta}. The group
$T^k$ acts freely on $\Phi^{-1}(0)$ and the toric manifold
corresponding to $\Delta$ is obtained as the symplectic quotient \(X
:= \C^N \mod T^k = \Phi^{-1}(0)/T^k.\) From the construction it may
be explicitly seen that $\Phi$ has a component which is proper and
bounded below, and the $T^k$-action on $\C^N$ has only finitely many
fixed point components. Hence we may apply
Theorem~\ref{theorem:kernelforT} to this setting.

Before proceeding to the computation of the kernel, we make a few
observations. First, there is a residual torus $T^n$ acting on $X$,
with moment map image $\Delta$. The facet
$F_i$ of $\Delta$ is exactly the $T^n$-moment map image of
\((\{z_i=0\} \cap \Phi^{-1}(0))/T^k.\) In particular, for a subset \(A
\subseteq \{1, 2, \ldots, N\},\) 
\begin{equation}\label{eq:capF}
\bigcap_{i \in A^c} F_i = \emptyset 
\quad\Longleftrightarrow\quad \C^A \cap \Phi^{-1}(0)
= \emptyset, 
\end{equation}
where $\C^A \subseteq \C^N$ denotes the coordinate subspace obtained
by setting $z_j = 0$ for $j \in A^c$. 
Second, we observe that 
$K^*_{T^k}(\C^N) \cong
R(T^k)$ may be expressed as the quotient of
$$K^*_{T^N}\left(\C^N\right) \cong
R\left(T^N\right) \cong \Z\left[x_1, \ldots, x_N, x_1^{-1}, \ldots, x_N^{-1}\right]$$
by the ideal 
\[
{\mathcal J} = \Bigl< x^{\alpha} -1  \Bigm\vert  \alpha \in \beta^*\bigl((\t^n)^*_{\Z}\bigr)
\subseteq (\t^N)^*_{\Z} \Bigr>,
\]
where $(\t^N)^*_{\Z} \cong \Z^N$ and $x^{\alpha} :=
x_1^{\alpha_1}x_2^{\alpha_2}\cdots x_N^{\alpha_N}$ for $x = (x_1,
\ldots, x_N),$ and $\alpha = (\alpha_1, \ldots, \alpha_N)$. 

We now proceed with the kernel computation. We first prove a technical
lemma which allows us to deal with the slight difficulty that an {\em
arbitrary} component of $\Phi$ need not be proper and bounded below. 

\begin{lemma}\label{lemma:gradient-flow}
Let $T^k$ be a compact torus acting linearly on $\C^N$ as
specified by~\eqref{eq:Delzant-Lie}, with moment map $\Phi: \C^N \to
(\t^k)^*$. Let $\Phi^{\xi}$ be a component of the moment map
and $c$ a critical value of $\Phi^{\xi}$. Let $\Sigma_c$ be the
corresponding critical set of $\Phi^{\xi}$ and let $D(\nu_c),
S(\nu_c)$ denote the disc and sphere bundles, respectively, of the negative normal
bundle $\nu_c$ to $\Sigma_c$ with respect to $\Phi^{\xi}$. Let $M^+_c
:= (\Phi^{\xi})^{-1}(-\infty, c+\varepsilon)$ and $M^-_c :=
(\Phi^{\xi})^{-1}(-\infty, c-\varepsilon)$. Then the
negative gradient flow corresponding to $\Phi^{\xi}$ with respect to
the standard K\"ahler metric on $\C^N$ gives a $T^k$-equivariant
retraction of the pair $(M_c^+, M_c^-)$ onto $(D(\nu_c), S(\nu_c))$.

\end{lemma}

\begin{proof}
  Let $F_{\Sigma}$ denote the union of the flow-up and flow-down sets
  contained in $M_c^+$ 
  associated to $\Sigma_c$. Since an 
arbitrary component $\Phi^{\xi}$ of the moment map may not be
  proper or bounded below on the noncompact space $\C^N$, we must show that any point \(z \in M_c^+
  \setminus F_{\Sigma}\) flows down to $M_c^-$ with respect to the
  negative gradient flow of $\Phi^{\xi}$. From this, the lemma
  follows. 

  In our case, the 
 vector field $- \grad(\Phi^{\xi})$ is complete by
  \cite[Proposition 2.1.21]{AbrMar78} so the flow
  is defined for all $t \in \R$. Suppose that \(z \in M_C^+ \setminus
  F_{\Sigma}\) and let $\gamma(t)$ for $t \in \R$ denote the negative gradient flow
  through $z$ with $\gamma(0) = z$. It will suffice to show that 
\[
\bigl|d\Phi^{\xi}_{\gamma(t)}(\gamma'(t))\bigr| =
\bigl|d\Phi^{\xi}_{\gamma(t)}(\grad \Phi^{\xi})_{\gamma(t)}\bigr| =
\bigl\|\grad(\Phi^{\xi})_{\gamma(t)}\bigr\|^2
\]
remains bounded away from $0$ as \(t \to \infty.\) On the other hand,
in the case of a linear action of $T^k$ on $\C^N$, we have 
\[
\bigl\|\grad(\Phi^{\xi})_{\gamma(t)}\bigr\|^2 = \bigl\|J\xi^{\sharp}_{\gamma(t)}\bigr\|^2
= \bigl\|\xi^{\sharp}_{\gamma(t)}\bigr\|^2 = \bigl\|\iota(\xi) \cdot \gamma(t)\bigr\|^2,
\]
where $J$ is the standard complex structure on $\C^N$.  

Given a $\xi \in \t^k$, let $\iota(\xi) = (\xi_1, \xi_2, \ldots,
\xi_N) \in \t^N \cong \R^N$ and let $A \subseteq \{1, 2, \ldots, N\}$ be the subset \(A =
\{j: \xi_j \neq 0\}.\) From the formula \(\grad(\Phi^{\xi})_y = J
\iota(\xi) y,\) for any \(y \in \C^N,\) we see that for $\ell \in A^c$, the
coordinate \(\gamma(t)_{\ell} = \gamma(0)_{\ell} = z_{\ell}\) is constant for all \(t \in
\R.\) Now let $C>0$ be such that \(|\iota(\xi)_j|^2 > C\) for all \(j
\in A.\) Let \(\pi_{\xi}: \C^N \to \C^A\) be the linear
projection which sends to $0$ all coordinates $z_{\ell}$ for \({\ell} \in A^c.\) 
Then we have the estimate 
\[
\bigl\|\grad(\Phi^{\xi})_{\gamma(t)}\bigr\|^2 > C\,\bigl\|\pi_{\xi}(\gamma(t))\bigr\|^2
\]
for all \(t >0\). Suppose now that 
\(\|\grad(\Phi^{\xi})_{\gamma(t)}\|^2 \to 0\) as \(t \to \infty.\)
Then \(\|\pi_{\xi}(\gamma(t))\|^2 \to 0\) also, and in particular
$\gamma(t)$ remains in a compact set. Hence $\gamma(t)$ converges. 
For a linear action of $T^k$ on $\C^N$, $\Phi^{\xi}$ has
only one critical component $\Sigma_c$, so $\gamma(t)$ must converge to
$\Sigma_c$, and we conclude $z \in F_{\Sigma}$. Therefore, if \(z \in
M_C^+ \setminus F_{\Sigma},\) then \(\|\grad(\Phi^{\xi})_{\gamma(t)}\|^2\)
is bounded away from $0$ for all $t$ and $z$ flows down to $M_c^-$. 
\end{proof}

With this lemma in hand, we may prove our main theorem.

\begin{theorem}\label{theorem:toric-variety}
Let $\Delta \subseteq \R^n$ be a compact Delzant polytope with facets
$\{F_i\}_{i=1}^N$ as specified in~\eqref{eq:def-Delta}.  Let \(X := \C^N \mod T^k\) be the symplectic
quotient of $\C^N$ corresponding to $\Delta$, as above. 
Then the ordinary $K$-theory ring $K^*(X)$ is given by
\[
K^*(X) \cong \Z[x_1, \ldots, x_N, x_1^{-1}, \ldots, x_N^{-1}]\Bigm/
{\mathcal I} + {\mathcal J}, 
\]
where 
\[
{\mathcal I} = \Bigl\langle {\prod}_{i \in S} (1 - x_i^{-1}) \Bigm\vert {\bigcap}_{i \in
    S} F_i = \emptyset \Bigr\rangle
\]
and 
\[
{\mathcal J} = \Bigl\langle x^{\alpha} -1 \Bigm\vert \alpha \in \beta^*((\t^n)^*_{\Z})
  \subseteq (\t^N)^*_{\Z} \cong \Z^N \Bigr\rangle.
\]

\end{theorem}

\begin{proof}
  We have already observed that $K^*_{T^N}(\C^N)/{\mathcal J}$ is isomorphic 
  to $K^*_{T^k}(\C^N)$, so to prove the theorem it suffices to
  describe $\ker(\kappa)$ in terms of the generators in
  $K^*_{T^N}(\C^N) \cong K^*_{T^N}(\pt)$. By
  Theorem~\ref{theorem:kernelforT}, we must compute $K_{\xi}$ for
  each $\xi \in Z$ as in~\eqref{eq:Z-def}. 

Since $T^k$ is acting linearly on $\C^N$, for each subset \(A
\subseteq \{1, 2, \ldots, N\}\) there is a unique critical value
$\xi_A \in \Phi(\C^A) = \Span_{\R_+}\{\alpha_i\}_{i \in A} - \eta$ 
associated to the linear subspace $\C^A$, and $Z$ consists of the
set of all $\xi_A$. (Here we denote by $\Span_{\R_+}\{\alpha_i\}_{i \in A}$ the set of linear
combinations of $\alpha_i$ with non-negative coefficients.) We note first that since $0$ is the only
$T^k$-fixed point of $\C^N$, and any class which restricts to $0$ at
the critical point $0$ is trivial, it suffices to consider $\xi_A$
such that $\langle \eta, \xi_A \rangle < 0$. 
It also suffices to consider only $\xi_A$ such that
\(\eta \not \in \Span_{\R_+} \{\alpha_i\}_{i \in A},\) which is
equivalent to 
\(\bigcap_{j \in A^c} F_j = \emptyset\) by~\eqref{eq:capF},
since otherwise \(\xi_A = 0\) and \(K_{\xi_A} = \{0\}.\) 

We now show that 
   for every \(\xi_A \in Z\) with \(\bigcap_{j \in A^c} F_j =
  \emptyset,\) the ideal
  $K_{\xi_A}$ is generated by an element of the form 
  \(\prod_{i \in S} (1 - x^{- 1}_i)\) for $S$ such that \(\bigcap_{i
   \in S} F_i = \emptyset.\) By definition of $Z$, the element 
  $\xi_A$ has the shortest
  norm in \(\Span_{\R_+}\{\alpha_i\}_{i \in A},\) so by construction $\xi_A$ has the
  property \(\langle \alpha_i, \xi_A \rangle = 0\) for all \(i \in A.\)
  Moreover, $\Phi$ is proper, so we have \(\langle \eta,
  \xi_A\rangle < 0.\) 
 By Theorem~\ref{theorem:injectivity}, in order to
  determine $K_{\xi_A}$ it suffices to examine the behavior of a class
   in $K_{\xi_A}$ at the
  unique $T^k$-fixed point $0 \in \C^N$. 
  Let
  $$M_0^+ := \bigl(\Phi^{\xi_A}\bigr)^{-1}(-\infty, +\varepsilon), \qquad
  M_0^- := \bigl(\Phi^{\xi_A}\bigr)^{-1}(-\infty, - \varepsilon)$$ for
  $0 < \varepsilon < \langle - \eta, \xi_A \rangle$. 
By Lemma~\ref{lemma:gradient-flow} we know that the pair 
  $(M_0^+, M_0^-)$ is $T^k$-equivariantly homotopic to 
 the disc and sphere bundle pair $(D(\nu_0),
  S(\nu_0))$. The negative normal bundle to the critical set of
  $\Phi^{\xi_A}$ is spanned by the coordinate
  directions $j$ such that \(\langle \alpha_j, \xi_A \rangle < 0.\)
  Let $S$ be the subset of such $j$. Then by construction we have 
  \(\eta \not \in \Span_{\R_+} \{\alpha_i\}_{i \in S^c}\) and thus \(\bigcap_{i \in S}
  F_i = \emptyset.\) Lemma~\ref{lemma:gradient-flow} and our choice of
  $\varepsilon$ also imply that 
  $M_0^-$ is $T^k$-equivariantly homotopic to $M_{\xi_A}$. These
  observations in addition to an analysis
  of the long exact sequence of the pair $(M_0^+, M_0^-)$
 imply that $K_{\xi_A}$ is
  generated by the class in
  $$K^*_{T^n}\left(\C^N\right) \cong K^*_{T^N}\left(\C^N\right)
  \bigm/ {\mathcal J} 
  \cong R\left(T^N\right)\bigm/{\mathcal J}$$ 
represented by $\prod_{i \in S} (1 - x_i^{- 1})$.

  Finally, we show that for every $S$ such that
  \(\bigcap_{i \in S} F_i = \emptyset,\) there exists $A$ such that 
  that the product $\prod_{i \in S} (1 - x^{- 1}_i)$ is contained in
  $K_{\xi_A}$. In fact it suffices to show this for minimal such $S$. 
By definition, 
for \(A := S^c\) we have that \(\eta \not \in \Span_{\R_+}
  \{\alpha_i\}_{i \in A}.\) Hence $\xi_A$ is non-zero and 
  $K_{\xi_A}$ is nontrivial. Moreover, since $S$ is minimal, we have 
 $\langle \xi_A, \alpha_i \rangle < 0$ for
  all \(i \in S.\) Hence each coordinate line for $i \in S$ is
  contained in the negative normal bundle to the critical set \(\C^S\)
  for the moment map component $\Phi^{\xi_A}$, and so by the argument
  above $K_{\xi_A}$ contains
  the product $\prod_{i \in S}(1 - x_i^{- 1})$. This completes the
  proof. 
\end{proof}

\def\cprime{$'$}

\end{document}